\ifx\shlhetal\undefinedcontrolsequence\let\shlhetal\relax\fi
\documentclass[11pt]{amsart}

\usepackage{amsmath}
\usepackage{amssymb}

\newtheorem{theorem}{Theorem}[section]
\newtheorem{claim}[theorem]{Claim}

\newtheorem{lemma}[theorem]{Lemma}

\newtheorem{corollary}[theorem]{Corollary}

\theoremstyle{definition}
\newtheorem{definition}[theorem]{Definition}

\theoremstyle{remark}
\newtheorem{remark}[theorem]{Remark}

\newcount\skewfactor
\def\mathunderaccent#1#2 {\let\theaccent#1\skewfactor#2
\mathpalette\putaccentunder}
\def\putaccentunder#1#2{\oalign{$#1#2$\crcr\hidewidth
\vbox to.2ex{\hbox{$#1\skew\skewfactor\theaccent{}$}\vss}\hidewidth}}

\def\smallbox#1{\leavevmode\thinspace\hbox{\vrule\vtop{\vbox
   {\hrule\kern1pt\hbox{\vphantom{\tt/}\thinspace{\tt#1}\thinspace}}
   \kern1pt\hrule}\vrule}\thinspace}

\newcommand{\then}{{\underline{then}}}
\newcommand{\Then}{{\underline{Then}}}

\def\qedref#1{$\qed_{\reforiginal{#1}}$}

\title{Separating properties for normal ultrafilters}
\author{Shimon Garti}
\address{Institute of Mathematics,
 The Hebrew University of Jerusalem,
 Jerusalem 91904, Israel}
\email{shimon.garty@mail.huji.ac.il}
\subjclass[2010]{Primary: 03E55, Secondary: 03E35, 03E05}
\keywords{Normal ultrafilters, measurable cardinals, compact cardinals, supercompact cardinals}

\begin{document}
\let\labeloriginal\label
\let\reforiginal\ref

\begin{abstract}
%%% Put abstract here:
Suppose $U$ is a normal ultrafilter on $\kappa$. We show that if $U$ contains the compact cardinals below $\kappa$ then every normal ultrafilter on $\kappa$ which contains the measurables, must contain the compact cardinals. On the other hand, if $U$ contains the supercompact cardinals below $\kappa$, then there is another normal ultrafilter $V$ which contains the compact but not the supercompact cardinals below $\kappa$. \newline
Supposons que $U$ est un ultrafiltre normal d\'efinie sur $\kappa$. Nous montrons que si $U$ contient l'ensemle des cardinaux compacts au-dessous de $\kappa$, chaque ultrafiltre normal sur $\kappa$ qui contient les cardinaux mesurables contient les cardinaux compacts. D'un autre c\^ot\'e, si $U$ contient les cardinaux supercompacts au-dessous de $\kappa$, il existe un autre ultrafiltre normal $V$ (sur $\kappa$) qui contient les cardinaux compacts mais ne contient pas les cardinaux supercompacts au-dessous de $\kappa$.
\end{abstract}

\maketitle

% start document here:
\newpage

\section{Introduction}

Suppose $\kappa$ is a measurable cardinal, and $U$ is a normal ultrafilter on $\kappa$. Let $p_0,p_1$ be two properties (defined by some formula) and set $A_0=\{\alpha<\kappa: p_0(\alpha)\},A_1=\{\alpha<\kappa: p_1(\alpha)\}$. We are interested in the possibility that $A_0\in U$ yet $A_1\notin U$.

Of course, this may happen trivially. For instance, let $p_0$ be the property of being weakly compact, and let $p_1$ be the opposite property of being not weakly compact. In this case, $A_0\in U$ and $A_1\notin U$ for every normal $U$ on every measurable $\kappa$.

To avoid this kind of trivialities, we shall ask for a measurable cardinal $\kappa$ which carries two normal ultrafilters $V,W$. For $V$ we shall have $A_0,A_1\in V$, and for $W$ we will be able to separate, namely $A_0\in W$ but $A_1\notin W$. The existence of $V$ shows that the separation is not construed on trivial grounds.

The following definition sets our terminology:

\begin{definition}
\label{cccognate} Cognate properties and separability. \newline
Assume $\kappa$ is a measurable cardinal, $p_0,p_1$ are two properties of ordinals, and $A_0=\{\alpha<\kappa: p_0(\alpha)\},A_1=\{\alpha<\kappa: p_1(\alpha)\}$.
\begin{enumerate}
\item [$(\aleph)$] $p_0$ and $p_1$ are cognates (for $\kappa$) if there exists a normal ultrafilter $U$ on $\kappa$ such that $A_0,A_1\in U$.
\item [$(\beth)$] $p_0$ and $p_1$ are separable (for $\kappa$) if they are cognates (for $\kappa$), but there is also a normal ultrafilter $W$ on $\kappa$ such that $A_0\in W$ and $A_1\notin W$.
\end{enumerate}
\end{definition}

We emphasize that the interesting problem is whether two properties are cognates and separable on the same $\kappa$. Moreover, we are looking for separability at every measurable cardinal $\kappa$. The following theorem of Solovay (see \cite{MR1994835}, 5.16) is a good example:

\begin{theorem}
\label{sssolovay} Normal ultrafilters which do not contain measurables. \newline
For every measurable cardinal $\kappa$ there exists a normal ultrafilter $U$ such that the set $\{\alpha<\kappa: \alpha$ is not a measurable cardinal $\}$ belongs to $U$.
\end{theorem}

\hfill \qedref{sssolovay}

In our terminology, we can define $p_0$ as being weakly compact, and $p_1$ as being measurable. The sets $A_0,A_1$ are defined respectively. The theorem of Solovay says that $p_0$ and $p_1$ are separable for every measurable cardinal $\kappa$ which carries a normal ultrafilter $U$ concentrating on the measurables.

Some point should be stretched here. If $\kappa$ is the first measurable cardinal, then every normal $U$ on $\kappa$ separates $p_0$ from $p_1$, as the set of measurables below $\kappa$ is empty. The theorem of Solovay begins to be interesting for higher measurables, on which a normal $V$ that contains the measurables exists as well. The same point will arise in the theorems below.

Let us fix three properties for the rest of the paper:
\begin{enumerate}
\item [$(\alpha)$] $p_0$ is the property of being measurable.
\item [$(\beta)$] $p_1$ is the property of being compact.
\item [$(\gamma)$] $p_2$ is the property of being supercompact.
\end{enumerate}

The corresponding sets $A_0,A_1,A_2$ are defined in accordance with $p_0,p_1,p_2$. We shall see that all of them can be cognates, but $p_0, p_1$ are never separable, and $p_1, p_2$ are separable. We indicate (in light of the above paragraph) that there exists a measurbale cardinal $\kappa$ and a normal $U$ on $\kappa$ such that $A_0\in U$ and $A_1\notin U$, but $p_0$ and $p_1$ are not cognates for this $\kappa$, and it will be the case for every normal ultrafilter on this $\kappa$. If there is one normal $U$ on some $\kappa$ which contains the compact cardinals, then for every normal ultrafilter on $\kappa$, we shall get $A_0\in U\Rightarrow A_1\in U$.

In fact, we shall prove the following:

\begin{theorem}
\label{mt1} measurable and compcat cardinals. \newline
Let $\kappa$ be a measurable cardinal, and let $U$ be a normal ultrafilter on $\kappa$ such that $A_0\in U$. \newline
\Then\ $A_1\in U$ iff $\kappa$ is a limit of compact cardinals.
\end{theorem}

\hfill \qedref{mt1}

This theorem means that the case of $A_0\in U$ and $A_1\notin U$ (for a normal $U$ on $\kappa$) happens only in the trivial case, in which the set of compact cardinals is bounded below $\kappa$. It follows that if one normal ultrafilter on $\kappa$ contains $A_1$, then all the normal ultrafilters on $\kappa$ contain $A_1$, provided that $A_0$ is in the ultrafilter, of course.

The other side of the coin is captured in the second theorem:

\begin{theorem}
\label{mt2} compact and supercompact cardinals. \newline
Let $\kappa$ be a measurable cardinal, and let $V$ be a normal ultrafilter on $\kappa$ such that $A_2\in V$. \newline
\Then\ there is another normal ultrafilter $W$ on $\kappa$ such that $A_1\in W$ but $A_2\notin W$.
\end{theorem}

\hfill \qedref{mt2}

In the second part of the paper, we apply separating properties to a problem concerning the number of normal ultrafilters which contain the non-measurables. We show that such properties enable us to create different normal ultrafilters by summing over a sequence of normal measures.

We conclude this introduction with an additional point of view. Magidor has proved in \cite{MR0429566} the consistency of the following:
\begin{enumerate}
\item [$(a)$] The first compact is also the first measurable.
\item [$(b)$] The first compact is also the first supercompact.
\end{enumerate}

In the phraseology of Magidor, these facts demonstrate an identity crisis of the compactness notion. We suggest the separability as another criterion for determining the identity of large cardinals. Due to the theorems above, compactness cannot be separated from measurability, but it is far from supercompactness (as much as normal ultrafilters are concerned).

We shall try to use standard notation, following \cite{MR1940513} and \cite{MR1994835}.
Recall that $\kappa$ is measurable if there exists a $\kappa$-complete non-principal ultrafilter on $\kappa$, so $\aleph_0$ is included, but throughout the paper we always mean that a measurable cardinal is uncountable. $\kappa$ is 1-extendible it there exists an elementary embedding $\jmath:V_{\kappa+1}\rightarrow V_{\jmath(\kappa)+1}$ so that $\kappa={\rm crit}(\jmath)$. $\kappa$ is compact if $\mathcal{L}_{\kappa\kappa}$ is compact (and here, again, we assume that $\kappa$ is uncountable). $\kappa$ is $\gamma$-supercompact if there exists an elementary embedding $\jmath:V\rightarrow M$ so that $\kappa={\rm crit}(\jmath)$ and $^\gamma M \subseteq M$. Finally, $\kappa$ is supercompact if $\kappa$ is $\gamma$-supercompact for every $\gamma$.

I thank professor Saharon Shelah for a very helpful discussion on the subject of this paper.

\newpage

\section{Separating theorems}

We commence with a general simple lemma:

\begin{lemma}
\label{llll} Suppose that:
\begin{enumerate}
\item [$(a)$] $\kappa$ is a measurable cardinal.
\item [$(b)$] $U$ is a normal ultrafilter on $\kappa$.
\item [$(c)$] $S\in U$.
\item [$(d)$] $T$ is an unbounded subset of $\kappa$.
\end{enumerate}
\Then\ the set $S'=\{\alpha\in S: \alpha$ is a limit of members from $T\}$ belongs to $U$.
\end{lemma}

\par\noindent\emph{Proof}.\newline
Assume toward contradiction that $S'\notin U$, so $S\setminus S'\in U$. For every $\alpha\in S\setminus S'$ we define $f(\alpha)={\rm sup}\{\beta<\alpha: \beta\in T\}$. By the definition of $S'$, each member $\alpha$ of $S\setminus S'$ is not a limit of members from $T$, hence $f(\alpha)<\alpha$. It follows that $f$ is a regressive function on $S\setminus S'$. As $U$ is normal (and $S\setminus S'\in U$) there is an ordinal $\gamma<\kappa$ and a subset $X\subseteq S\setminus S', X\in U$ such that:

$$
\alpha\in X \Rightarrow f(\alpha)=\gamma
$$

But $X\in U$, so in particular $X$ is unbounded in $\kappa$. It follows that $T\subseteq \gamma+1$, in contrary to assumption $(d)$ of the lemma.

\hfill \qedref{llll}

Menas proved (in \cite{MR0357121}) that if $\kappa$ is a measurable cardinal, limit of compact cardinals, then $\kappa$ is a compact cardinal. Combining this theorem with the lemma above, we can conclude:

\begin{claim}
\label{ccccom} Suppose that:
\begin{enumerate}
\item [$(a)$] $\kappa$ is a measurable cardinal.
\item [$(b)$] $U$ is a normal ultrafilter on $\kappa$.
\item [$(c)$] $U$ contains the set of measurables below $\kappa$.
\end{enumerate}
\Then\ $U$ contains the set of compact cardinals below $\kappa$ iff $\kappa$ is a limit of compact cardinals.
\end{claim}

\par\noindent\emph{Proof}. \newline
Clearly, if the set of compact cardinals below $\kappa$ is bounded then it is not in $U$. Assume that this set is unbounded in $\kappa$. Denote the set of measurable cardinals below $\kappa$ by $S$ (so $S\in U$ by assumption $(c)$), and the set of compact cardinals below $\kappa$ by $T$.

Applying Lemma \ref{llll} we know that the set $S'=\{\alpha\in S: \alpha$ is a limit of compact cardinals $\}$ belongs to $U$. By the theorem of Menas, $S'$ is a set of compact cardinals, hence $S'\subseteq T$. It follows that $T\in U$, as required.

\hfill \qedref{ccccom}

The meaning of this claim is that the properties $p_0,p_1$ are not separable. In other words, if $\kappa$ carries a mormal ultrafilter $U$ which contains the compact cardinlas below $\kappa$, then every normal ultrafilter $V$ which contains the measurables must agree with $U$ about the compact cardinlas below $\kappa$.

\begin{remark}
\label{mmmm} It is consistent that $\kappa$ is a measurable cardinal, $U$ is a normal ultrafilter on $\kappa$, the set of measurables belongs to $U$, yet $U$ does not contain the set of compact cardinlas below $\kappa$.

Indeed, take the model of Magidor, in which the first compact $\kappa$ is also the first supercompact. Being supercompact entails the existence of a normal $U$ which contains the measurables, and being the first compact means that the set of compact cardinlas below $\kappa$ is empty. Nevertheless, in this case every normal $V$ on $\kappa$ which contains the measurables, behaves like $U$.
\end{remark}

\hfill \qedref{mmmm}

Our goal is to prove that $p_1$ and $p_2$ are separable for every relevant $\kappa$. We shall prove this by induction on the measurables, and it holds for every measurable cardinal $\kappa$ for which $p_1$ and $p_2$ are cognates. Let us phrase the accurate assertion:

\begin{theorem}
\label{mt} Compactness and Supercompactness. \newline
If $\kappa$ is a measurable cardinal, $V$ is a normal ultrafilter on $\kappa$, and $A_2\in V$, \then\ there is a normal ultrafilter $W$ on $\kappa$ such that $A_1\in W$ but $A_2\notin W$.
\end{theorem}

\par \noindent \emph{Proof}. \newline
By induction on the measurables which carry a normal ultrafilter $V$ so that $A_2\in V$. Suppose we have arrived at $\kappa$ which satisfies this assumption, and the theorem holds below $\kappa$. For every measurable cardinal $\alpha$ below $\kappa$ we choose a normal ultrafilter $U_\alpha$ on $\alpha$ as follows:

If there is no normal ultrafilter on $\alpha$ which contains the measurables, then let $U_\alpha$ be any normal measure on $\alpha$. If there is a normal ultrafilter on $\alpha$ which contains the measurables, we choose $U_\alpha$ such that $A_0\in U_\alpha$ and $A_2\notin U_\alpha$. This is possible either by the induction hypothesis or if there is no normal ultrafilter on $\alpha$ such that $A_2\in U_\alpha$ at all.

Notice that $\kappa$ is a limit of compact cardinals, by the assumptions of the theorem.
Let $V$ be any normal measure on $\kappa$ such that $A_0\in V$. If $A_2 \notin V$ then we are done (as $A_1\in V$ by virtue of Claim \ref{ccccom}), so assume $A_2\in V$. We shall elicit a new normal ultrafilter $W$ on $\kappa$, as required. For every subset $X\subseteq \kappa$ we define:

$$
X\in W \Leftrightarrow \{\alpha<\kappa: X\cap\alpha\in U_\alpha\}\in V.
$$

It is easily verified that $W$ is an ultrafilter on $\kappa$. Let us show that $W$ is normal. Suppose $X\in W$ and $f$ is a regressive function on $X$. Denote the set $\{\alpha<\kappa:X\cap\alpha\in U_\alpha\}$ by $S_X$. It follows from the definition of $W$ that $S_X\in V$. As $f\upharpoonright(X\cap\alpha)$ is regressive for every $\alpha\in S_X$, there exist a set $Y_\alpha\subseteq X\cap\alpha, Y_\alpha\in U_\alpha$ and an ordinal $\gamma_\alpha<\alpha$ so that $f\upharpoonright Y_\alpha=\{\gamma_\alpha\}$. The mapping $\alpha\mapsto \gamma_\alpha$ is a regressive function on $S_X$, so (by the normality of $V$) there is a set $S\subseteq S_X, S\in V$ and a fixed ordianl $\gamma<\kappa$ such that $\alpha\in S\Rightarrow\gamma_\alpha=\gamma$. Let $Y$ be $\bigcup\{Y_\alpha: \alpha\in S\}$. By the considerations above, $Y\in W$ (as $S_Y\supseteq S\in V$) and $f\upharpoonright Y=\{\gamma\}$. It follows that $W$ is normal, as required.

We claim that $W$ satisfies the theorem. Indeed, $A_2\cap\alpha\notin U_\alpha$ for every measurable $\alpha$ in $A_0$ (by the choice of $U_\alpha$). As $A_0\in V$ it means that $A_2\notin W$.

On the other hand, $A_0\cap\alpha\in U_\alpha$ for every supercompact $\alpha$ in $A_2$ (as every supercompact $\alpha$ carries a normal ultrafilter which concentrates on the measurables), hence $A_0\in W$ (recall that $A_2\in V$). As $\kappa$ is a limit of compact cardinals (recall that $A_2\in V$ and every supercompact is compact) we conclude from Claim \ref{ccccom} that $A_1\in W$ as well, so we are done.

\hfill \qedref{mt}

One may wonder if the theorem of Menas can be rephrased for measurables and supercompact cardinals or even for compact and supercompact cardinals. We can infer from the theorem above that if $\kappa$ is measurable and a limit of supercompact cardinals, then $\kappa$ need not be a supercompact cardinal (but notice that $\kappa$ is a compact cardinal). In fact, the amount of measurables, limit of supercompact cardinals, which are not supercompact by themselves, is large in the following sense:

\begin{corollary}
\label{nnnnmenas} Assume that:
\begin{enumerate}
\item [$(a)$] $\kappa$ is a measurable cardinal.
\item [$(b)$] $\kappa$ is a limit of supercompact cardinals.
\item [$(c)$] $U$ is a normal ultrafilter on $\kappa$.
\item [$(d)$] $A_0\in U$.
\end{enumerate}
\Then\ there exists a normal ultrafilter $V$ on $\kappa$ such that the set of measurables limit of supercompact cardinals which are not supercompact, belongs to $V$. Moreover, the members of this set are compact cardinals.
\end{corollary}

\par\noindent\emph{Proof}.\newline
By Theorem \ref{mt} we can create $V$ such that $A_0\in V$ and $A_2\notin V$. Let $A'\subseteq A_0$ be the set of measurables which are limit of supercompact cardinals. By Lemma \ref{llll} we know that $A'\in V$ (assumption $(b)$ is used here, and $A_2$ plays the role of $T$ in the lemma).

As $A_2\notin V$ we know that $A'\setminus(A'\cap A_2)\in V$ Notice that the members of this set are not supercompact although each of them is a limit of supercompact cardinals. By the theorem of Menas, the members of $A'\setminus(A'\cap A_2)$ are compact cardinals. Hence the set $A'\setminus(A'\cap A_2)$ completes the proof of the corollary.

\hfill \qedref{nnnnmenas}

However, if there is a normal ultrafilter $U$ on $\kappa$ so that $A_2\in U$, then the set of supercompact cardinals which are limit of supercompact cardinals belongs to $U$ (by Lemma \ref{llll}). The existence of $U$ is coherent, of course, with Corollary \ref{nnnnmenas}.
It should be observed that under stronger assumption we can get the following parallel to the theorem of Menas:

\begin{claim}
\label{mmmmenas} Suppose that:
\begin{enumerate}
\item [$(a)$] $\kappa$ is a measurable cardinal.
\item [$(b)$] $U$ is a normal ultrafilter on $\kappa$.
\item [$(c)$] $A_2\in U$.
\end{enumerate}
\Then\ $\kappa$ is supercompact.
\end{claim}

\par\noindent\emph{Proof}.\newline
Let $\gamma$ be any ordinal above $\kappa$. We have to show that $\kappa$ is $\gamma$-supercompact, i.e., that there is a fine and normal measure on $\mathcal{P}_\kappa\gamma$. For every $\alpha\in A_2$ we choose a fine and normal ultrafilter $U_\alpha$ on $\mathcal{P}_\alpha\gamma$. Now, for every $x\subseteq \mathcal{P}_\kappa\gamma$ we define:

$$
x\in V \Leftrightarrow \{\alpha<\kappa: x\cap\mathcal{P}_\alpha\gamma\in U_\alpha\}\in U.
$$

We claim that $V$ is a fine and normal ultrafilter on $\mathcal{P}_\kappa\gamma$. The fact that $V$ is an ultrafilter is routine. We know that $V$ is normal since it results as a sum of normal ultrafilters over some normal ultrafilter. Let us show that $V$ is fine.

Fix any $i\in\gamma$ and let $S_i$ be the set $\{x\in \mathcal{P}_\kappa\gamma: i\in x\}$. We must prove that $S_i\in V$. Observe that $S_i\cap \mathcal{P}_\alpha\gamma = \{x\in \mathcal{P}_\kappa\gamma: |x|<\alpha\wedge i\in x\}=\{x\in \mathcal{P}_\alpha\gamma: i\in x\}\in U_\alpha$.
It follows that $S_i\cap \mathcal{P}_\alpha\gamma\in U_\alpha$ for every $\alpha\in A_2$. This is suffice, as $A_2\in U$.

\hfill \qedref{mmmmenas}

The main theorem above can be rephrased in a more general setting. We shall prove, on the basis of this generalization, that some higher properties are separable:

\begin{theorem}
\label{gggg} The generalized separation. \newline
Let $\kappa$ be a measurable cardinal, and let $p,r$ be two properties such that $A_r\subseteq A_p$, and there is a normal ultrafilter $V$ on $\kappa$ such that $A_r\in V$. Assume that $r(\alpha)\Rightarrow$ there is a normal $U_\alpha$ on $\alpha$ so that $A_p\in U_\alpha$. \newline
\Then\ the properties $p$ and $r$ are separable.

In particular, the properties of measurability and $1$-extendibility are separable, as well as the properties $1$-extendibility and $2^\kappa$-supercompactness.
\end{theorem}

\par\noindent\emph{Proof}.\newline
Let $X$ be the collection of measurables which carry a normal ultrafilter $V$ so that $A_r\in V$. We procceed by induction on the members of $X$. Assume $\kappa\in X$, and the theorem holds on the members of $X$ below $\kappa$. Fix $V$ on $\kappa$ so that $A_r\in V$. We shall define $W$, another normal ultrafilter on $\kappa$, such that $A_p\in W$ yet $A_r\notin W$.

For this end, we choose a sequence of normal ultrafilters $\langle U_\alpha: \alpha<\kappa\rangle$ (for every measurable $\alpha$) as follows. If there is no normal ultrafilter on $\alpha$ which contains $A_p$ then let $U_\alpha$ be any normal ultrafilter. Otherwise, let $U_\alpha$ satisfy $A_p\in U_\alpha \wedge A_r\notin U_\alpha$. This is possible by the induction hypothesis.

For every $y \subseteq \kappa$ define:

$$
y\in W \Leftrightarrow \{\alpha<\kappa: y\cap\alpha\in U_\alpha\}\in V.
$$

As in the main theorem above, $W$ is normal, and $A_p\in W \wedge A_r\notin W$, so we are done.

Recall that if $\kappa$ is $2^\kappa$-supercompact then $\kappa$ carries a normal ultrafilter $U$ which contains the set of $1$-extendible cardinals below $\kappa$. Inasmuch as every $1$-extendible cardinal carries a normal ultrafilter which contains the measurables, we know that measurability and $1$-extendibility are separable. Likewise, measurability and $2^\kappa$-supercompactness are separable. Hence, one can prove a parallel to Theorem \ref{mt} with respect to measurability and $2^\kappa$-supercompactness.

The same holds for higher notions of large cardinals. For instance, each extendible cardinal carries a normal ultrafilters which contains the supercompact cardinals. It follows that extendibility and supercompactness are separable.

\hfill \qedref{gggg}

The theorems of this section provide some stable distinction between compact and supercompact cardinals. This is natural when the defining criterion is related to normality of ultrafilters. Indeed, the notion of normality plays a central role in the definition of supercompactness, but not in the definitions of measurability or compactness.

\newpage

\section{Normal measures containing the non-measurables}

As mentioned in the previous section, Solovay proved that each measurable cardinal $\kappa$ carries a normal measure which contains the non-measurables below $\kappa$. Clearly, there is no parallel theorem for normal measures which concentrate on the measurables. If $\kappa$ is the first measurable (or any measurable cardinal which is not a limit of measurable cardinals) then the set of measurables below $\kappa$ cannot be a member of a normal ultrafilter.

Nonetheless, if we assume more than measurability then we have normal ultrafilters containing the measurables. It suffices that $\kappa$ is $2^\kappa$-supercompact to prove the existence of such an ultrafilter. Moreover, there are many distinct normal ultrafilters on $\kappa$ which concentrate on the measurables if we assume enough supercompactness.

Despite that fact, we do not know if there are two distinct normal ultrafilters on $\kappa$ which contain the non-measurables. This is labeled as an open problem in \cite{MR1994835}. It seems that separability properties below the notion of measurability cannot be established, so a direct approach for creating two distinct normal ultrafilters fails. The essential obstacle is that weak properties (namely, weaker than measurability) can hardly be separated in the sense of Theorem \ref{gggg}. We cannot assume that $r(\alpha)\Rightarrow$ there is a normal $U_\alpha$ on $\alpha$ so that $A_p\in U_\alpha$, since there is no normal ultrafilter on $\alpha$ at all if the property $r$ is weaker than measurability.

Moreover, it may happen that the normal ultrafilter ensured by Theorem \ref{sssolovay} of Solovay is the only normal ultrafilter on $\kappa$ (see \cite{MR0277346}). If so, there are no two distinct normal ultrafilters (on $\kappa$) containing the non-measurables.
Nevertheless, separability can be invoked in a slight different manner.
In order to use separability for producing two distinct normal ultrafilters which do not contain the measurable cardinals, we shall act as follows.

In the theorems below we shall see that separating properties \emph{above measurability} enable us to create distinct normal ultrafilters which contain the non-measurables (but we do not resolve the above open problem). Let $\kappa$ be a measurbale cardinal, and for every measurable $\alpha<\kappa$ choose a normal ultrafilter $U_\alpha$ which does not contain the measurables.
We shall try to define sums of the same sequence $\langle U_\alpha: \alpha<\kappa\rangle$ over distinct normal ultrafilters on $\kappa$ which contain the measurables. Unfortunately, there is no reason that this process will be one-to-one, but separating properties may help:

\begin{theorem}
\label{aaaa} Ultrafilters which do not contain the measurables. \newline
Let $\kappa$ be a measurable cardinal, and let $p,q$ be two properties that are separable on $\kappa$. Let $U_\alpha$ be a normal ultrafilter on $\alpha$ which does not contain the measurables, for every measurable $\alpha<\kappa$, and let $S_X$ be $\{\alpha<\kappa: X\cap \alpha\in U_\alpha\}$ for every $X\subseteq \kappa$. \newline
If there is $X\subseteq \kappa$ so that $S_X=A_p\setminus A_q$, then there are (at least) two distinct normal ultrafilters on $\kappa$ which do not contain the measurables.
\end{theorem}

\par \noindent \emph{Proof}. \newline
Let $\{V_\gamma: \gamma<\gamma^*\}$ be the set of all normal ultrafilters on $\kappa$ which contain the measurables. Let $D$ be $\bigcap \{V_\gamma: \gamma<\gamma^*\}$. Observe that $D$ is a normal filter (as an intersection of normal filters) but not an ultrafilter (unless $\gamma^*=1$). Denote the complement ideal by $I$.

We claim that if there is a subset $X\subseteq \kappa$ such that $S_X\notin D\cup I$ then we are done. Indeed, if such an $X$ exists, then we may find two distinct ultrafilters $V_0,V_1$ so that $S_X\in V_0$ and $S_X\notin V_1$. We define $W_0,W_1$ as usual, namely $X\in W_\ell \Leftrightarrow S_X\in V_\ell$ for $\ell\in\{0,1\}$. It follows that each $W_\ell$ is a normal ultrafilter which do not contain the measurables (as every $U_\alpha$ runs away from the measurables, and the set of $\alpha$-s belongs to $V_\ell$).
Moreover, $W_0\neq W_1$, as demonstrated by $X$. Indeed, $X\in W_0$ as $S_X\in V_0$, but $X\notin W_1$ as $S_X\notin V_1$.

Now choose $X\subseteq\kappa$ so that $S_X=A_p\setminus A_q$. Let $V_0,V_1$ be normal ultrafilters on $\kappa$ such that $A_p\setminus A_q\in V_0$ and $A_p,A_q\in V_1$. This can be done by the assumptions of the Theorem. It follows that $S_X\notin V_1$ (as $A_q\in V_1$) and $S_X\in V_0$, so we are done.

\hfill \qedref{aaaa}

\begin{remark}
\label{oooo} The opposite direction holds as well. If $\langle U_\alpha:\alpha<\kappa\rangle$ is fixed, $V_0\neq V_1$ are normal ultrafilters on $\kappa$ which contain the measurables, $W_0,W_1$ are the sums of $\langle U_\alpha:\alpha<\kappa\rangle$ over $V_0,V_1$ and $W_0\neq W_1$, \then\ there exists $X\subseteq\kappa$ such that $S_X,\kappa\setminus S_X\notin V_0\cap V_1$.
\end{remark}

\hfill \qedref{oooo}

Based on Theorem \ref{aaaa}, we shall introduce a theorem about normal ultrafilters which do not contain the measurables. The idea is to use the definition of Anti-Limit-Condition (ALC, to be defined below). The strong version of this notion yields distinct ultrafilters in ZFC, and we indicate that the weaker version enables us to force the existence of distinct normal measures on $\kappa$ (which do not contain the measurables) whenever $\kappa$ is supercompact. The method is to procure a subset $X\subseteq\kappa$ so that $S_X$ belongs to one normal ultrafilter on $\kappa$ but not to some other such ultrafilter.

\begin{definition}
\label{aallcc} The Anti-Limit-Condition.
\begin{enumerate}
\item [$(a)$] Let $\alpha$ be a measurable cardinal, limit of measurable cardinals. Let $\langle U_\beta:\beta\leq\alpha\wedge\beta$ is a measurable cardinal$\rangle$ be a sequence of normal ultrafilters, each $U_\beta$ is defined on the measurable cardinal $\beta$. \newline
We say that $\langle U_\beta:\beta\leq\alpha\rangle$ has the ALC if one can find a sequence of sets $\langle X_\beta:\beta<\alpha\rangle$ such that $X_\beta\in U_\beta$ for every $\beta<\alpha$ but $\bigcup\{X_\beta:\beta<\alpha\}\notin U_\alpha$.
\item [$(b)$] Let $\kappa$ be a measurable cardinal, $T_1$ an unbounded subset of $\kappa$ which is a set of measurable cardinals and $T=\{\alpha:\alpha$ is a measurable cardinal limit of members from $T_1\}$. \newline
We say that $\kappa$ satisfies the ALC with respect to $T_1$ if there is a sequence $\langle U_\beta:\beta\in T_1\cup T\rangle$ of normal ultrafilters (each $U_\beta$ is defined on the measurable cardinal $\beta$) such that the sequence $\langle U_\beta:\beta\leq\alpha\rangle$ has the ALC whenever $\alpha\in T$.
\item [$(c)$] Let $\kappa$ be a measurable cardinal, $T_1$ an unbounded subset of $\kappa$ which is a set of measurable cardinals and $T=\{\alpha:\alpha$ is a measurable cardinal limit of members from $T_1\}$. \newline
We say that $\kappa$ satisfies the uniform ALC with respect to $T_1$ if it satisfies the ALC with respect to $T_1$ (as exemplified by $\langle U_\beta:\beta\in T_1\cup T\rangle$) and one can choose $\langle X_\beta:\beta<\kappa\rangle$ such that $X_\beta\in U_\beta$ for every $\beta\in T_1\cup T$, and $\bigcup\{X_\beta: \beta<\alpha\}\notin U_\alpha$ whenever $\alpha\in T$.
\end{enumerate}
\end{definition}

\begin{remark}
\label{alcccc} Part $(a)$ defines a property of the sequence of ultrafilters. Part $(b)$ deals with a property of the measurable $\kappa$ itself, demanding that the property of part $(a)$ holds for every measurable cardinal limit of measurable cardinals (in the set $T$) below $\kappa$. The extra requirement in part $(c)$ is that the sets in the ultrafilters are taken from some fixed sequence of sets, while in part $(b)$ we may change them for every specific $\alpha$.
\end{remark}

\hfill \qedref{alcccc}

The uniform ALC property (of part $(c)$) provides the tools for the construction of two distinct ultrafilters which do not contain the measurables in ZFC, as follows:

\begin{theorem}
\label{uuuu} Uniform ALC and distinct normal measures. \newline
Let $\kappa$ be a measurable cardinal, satisfies the uniform ALC. \newline
Suppose $T_0,T_1$ are separable for $\kappa$, $T_0\subseteq T_1$ and both consist of measurable cardinals. \newline
\Then\ one can elicit two distinct normal ultrafilters on $\kappa$ which do not contain the measurable cardinlas below $\kappa$, by summing over a fixed sequence of ultrafilters.
\end{theorem}

\par\noindent\emph{Proof}.\newline
Set $T=\{\alpha: \alpha$ is a measurable cardinal, limit of members from $T_1\}$. Let $\langle U_\beta:\beta\in T_1\cup T\rangle$ exemplify the uniform ALC on $\kappa$ with respect to $T_1$. We compute $S_X$ according to this sequence, i.e., $S_X=\{\beta<\kappa:X\cap\beta\in U_\beta\}$ for every $X\subseteq\kappa$.
Choose a sequence $\langle X_\beta: \beta\in T_1\cup T\rangle$ as in Definition \ref{aallcc}(c). Define:

$$
X=\bigcup\{X_\beta: \beta\in T_1\cup T\}
$$

We claim that $S_X=T_1$. For this, consider an ordinal $\beta\in T_1$. If follows (from the very definition of the set $X$) that $X_\beta\subseteq X\cap\beta$. As $X_\beta\in U_\beta$ we have $X\cap\beta\in U_\beta$, hence $\beta\in S_X$. On the other hand, if $\beta\notin T_1$ then we may assume (without loss of generality) that $\beta\in T$ (notice that we ignore measurable cardinals outside $T_1\cup T$, as this set is negligible). By the ALC we know that $X\cap\beta=\bigcup\{X_\gamma:\gamma<\beta\}\notin U_\beta$, so $\beta\notin S_X$. Now by virtue of Theorem \ref{aaaa} we are done.

\hfill \qedref{uuuu}

Essentially, the same idea may be rendered in order to force the existence of distinct normal ultrafilters which do not contain the measurables. We have to assume merely ALC (i.e., part $(b)$ of Definition \ref{aallcc}). Suppose $\kappa$ is a supercompact cardinal. By the preparatory forcing of Laver (see \cite{MR0472529}) we may assume that $\kappa$ is indestructible.

We begin with two distinct normal ultrafilters $V_0,V_1$ which contain the measurables, and separate the properties $T_0,T_1$. Now we can employ the so-called generalized Mathias forcing (see \cite{MR2927607}, for an explicit description of it) to force the existence of a pair $(s,x)$ so that $s\subseteq S_x\subseteq T_1$. If we use $V_0,V_1$ for summing over a fixed sequence of normal ultrafilters which satisfies \ref{aallcc}(b), we produce $W_0,W_1$ as required (i.e., normal ultrafilters which do not contain the measurables such that $x\in W_1$ yet $x\notin W_0$).

Anyway, as Saharon Shelah indicated, much simpler arguments may be invoked in order to force distinct normal ultrafilters which do not contain the measurables. It means that the original problem is interesting only in the frame of ZFC.

\begin{remark}
\label{ssss}
The existence of two distinct normal ultrafilters which do not contain the measurables need not come from summing over normal ultrafilters which concentrate on the measurables. By Kunen-paris (see \cite{MR0277381}) one can force the existence of many normal ultrafilters on the first measurable cardinal. All of them do not contain the measurabels, but none of them comes from the process described in the theorems of this paper.
\end{remark}

\hfill \qedref{ssss}

\newpage

\bibliographystyle{amsplain}
\bibliography{arlist}
%liste

\end{document}